\input amstex
\documentstyle{amsppt}
\magnification=\magstep1 \NoRunningHeads
\topmatter

\title
 Mixing actions of $0$-entropy for countable amenable groups
 \endtitle

\author
Alexandre I. Danilenko
\endauthor

\email
alexandre.danilenko@gmail.com
\endemail

\address
 Institute for Low Temperature Physics
\& Engineering of National Academy of Sciences of Ukraine, 47 Lenin Ave.,
 Kharkov, 61164, UKRAINE
\endaddress
\email alexandre.danilenko\@gmail.com
\endemail

\abstract
It is shown that each discrete countable infinite amenable group admits a 0-entropy mixing 
action on a standard probability space.
\endabstract

 \loadbold

\endtopmatter

\NoBlackBoxes

\document

\head 0. Introduction
\endhead

Let $G$ be an amenable  discrete  infinite  countable  group.
We recall that a measure preserving action $T=(T_g)_{g\in G}$ of $G$ on a standard probability space  $(X,\goth B,\mu)$ is called {\it mixing}
if $\mu(T_gA\cap B)\to\mu(A)\mu(B)$ as $g\to\infty$ for all measurable subsets $A,B\subset X$.
The Kolmogoroff-Sinai entropy for measure-preserving transformations (i.e. $\Bbb Z$-actions) was extended to the actions of locally compact amenable groups in \cite{OrWe}.
Since $0$-entropy actions are considered as ``deterministic'' while mixing actions are considered as ``chaotic'', it is natural to ask are there actions which possess the two properties   simultaneously?
Such examples for $\Bbb Z$-actions were found  first in the class of Gaussian transformations (see \cite{New}) and later in the class of rank-one transformations (see \cite{Or}, \cite{Ad}, \cite{CrSi} and references therein).
The rank-one analogues of the later family  were constructed   
for $G=\Bbb R^{d_1}\times\Bbb Z^{d_2}$ in \cite{DaSi},  for $G$ being a countable direct sum of finite groups in \cite{Da2} and, more generally,   for $G$ being a locally normal countable group
 in  \cite{Da3}.
Dan Rudolph asked\footnote{At an AMS meeting in early 90's.} whether {\it each} amenable countable group $G$ has a mixing   action of zero entropy? 
The purpose of this work is to answer this question affirmatively.

\proclaim{Theorem 0.1} There is a 0-entropy mixing free\footnote{We stick to the free actions to avoid degeneracy and triviality like the following:  if $(T_g)_{g\in G}$ is a 0-entropy mixing action of $G$,  $F$ is a finite group  and  $\widetilde T_{g,f}:=T_g$ for $(g,f)\in G\times F$ then  $(\widetilde T_{g,f})_{(g,f)\in G\times F}$ is  a 0-entropy mixing  action of $G\times F$. We consider this non-free action as   degenerated and non-interesting.} probability preserving action of $G$.
\endproclaim

Since the proof of the theorem is based on the Poisson suspensions of infinite measure preserving actions we introduce the necessary   definitions in the next  section (see \cite{Ne}, \cite{Roy}, \cite{Ja--Ru} for more detail).
We note that the mixing property follows essentially from the construction elaborated in \cite{Da4}.
Thus we have only to show here how to modify that  construction to achieve the freeness\footnote{The freeness of the $H_3(\Bbb R)$-actions considered in \cite{Da4} follows from some properties of  group actions that are specific to actions of  connected nilpotent Lie groups.} and 0-entropy of the action in question.

\head 1. Poisson suspensions\endhead

 Let $X$  be a locally compact noncompact Cantor space (i.e. $0$-dimensional, without isolated points).
Denote by $C_{00}(X)$ the  vector space of real valued functions  on $X$ with compact support.
This space is endowed with the usual locally convex topology, i.e. the topology of uniform convergence on the compact subsets of $X$.
The dual $C_{00}(X)'$ is called the space of (real) Radon measures on $X$.
We are interested in the cone $X^*\subset C_{00}(X)'$ 
 of nonnegative Radon measures on $X$.
Furnish $X^*$ with the Borel $\sigma$-algebra $\goth B^*$ generated by the $*$-weak topology related to the duality $\langle C_{00}(X),C_{00}(X)'\rangle$.
We note that $\goth B^*$ is also the   Borel $\sigma$-algebra generated by the strong topology related to this duality.
 Since the strong topology is Polish and $X^*$ is closed in $C_{00}(X)'$, it follows that
  $(X^*,\goth B^*)$ is a standard Borel space.

Denote by  $\Cal K$ the set of all compact open subsets of $X$.
 Of course, $\Cal K$ is infinite but countable.
We also note that  $\goth B^*$ is the smallest $\sigma$-algebra on $X^*$
such that for each $K\in\Cal K$, the mapping $N_K: X^*\ni x^*\mapsto x^*(K)\in\Bbb R_+$ is measurable.
Let $\mu^*$ be the only measure on $(X^*,\goth B^*)$ such that
\roster
\item"---"
$\mu^*\circ N_K^{-1}$ is the Poisson distribution with parameter $\mu(K)$ and
\item"---" the random variables $N_{K_1},\dots, N_{K_n}$  on $(X^*,\mu^*)$ are independent for each countable collection of mutually disjoint subsets $K_1,\dots,K_m\in\Cal K$.
\endroster
Then $(X^*,\goth B^*,\mu^*)$ is a standard probability space.
To define $\mu^*$ rigorously we denote by $\goth F$ the set of all finite collections of mutually disjoint nonempty compact open subsets of $X$.
For $\Cal F_1,\Cal F_2\in\goth F$, we write $\Cal F_2\succ\Cal F_1$ if each element of $\Cal F_1$ is a union of some elements from $\Cal F_2$.
Then $(\goth F,\succ)$ is a directed partially ordered set.
Given $K\in\Cal K$, we define a probability measure $\nu_K$ on $\Bbb R_+$ by setting 
$$
\nu_K(A):=\sum_{i\in A\cap\Bbb Z_+}\frac{e^{-\mu(K)}\mu(K)^i}{i!}
$$
for each Borel subset $A\subset\Bbb R_+$.
We see, in particular, that $\nu_K$ is supported on $\Bbb Z_+$.
Since $\Bbb R_+$ is an additive semigroup, the convolution of probability measures on $\Bbb R^*$ is well defined. 
It is easy to verify that if $K_1,K_2\in\Cal K$ and $K_1\cap K_2=\emptyset$ then $\nu_{K_1}*\nu_{K_2}=\nu_{K_1\sqcup K_2}$.
Given $\Cal F\in\goth F$, we denote by $\Bbb R_+^{\Cal F}$ the set of all mappings $x:\Cal F\ni K\mapsto x(K)\in\Bbb R_+$.
We define a measure  $\nu^{\Cal F}$ on $\Bbb R_+^{\Cal F}$ as the direct product   $\nu^{\Cal F}:=\bigotimes_{K\in\Cal F}\nu_K$.
If $\Cal F_1,\Cal F_2\in\goth F$ and $\Cal F_2\succ\Cal F_1$
we define a mapping $\pi_{\Cal F_1}^{\Cal F_2}:\Bbb R_+^{\Cal F_2}\to\Bbb R_+^{\Cal F_1}$ by setting
$$
\bigg(\pi_{\Cal F_1}^{\Cal F_2}(x)\bigg)(K)=\left(\sum_{\{R\in\Cal F_2\mid R\subset K\}} x(R)\right)\quad\text{for each $K\in\Cal F_1$ and } x\in\Bbb R_+^{\Cal F_2}.
$$
The aforementioned convolution property of $\nu_K$ implies that  $\nu^{\Cal F_2}\circ(\pi_{\Cal F_1}^{\Cal F_2})^{-1}=\nu^{\Cal F_1}$.
Thus $\{(\Bbb R^{\Cal F_1},\nu^{\Cal F_2},\pi_{\Cal F_1}^{\Cal F_2})_{\Cal F_2\succ\Cal F_1}\mid \Cal F_1,\Cal F_2\in\goth F\}$ is a projective system of probability spaces.
Hence the projective limit $(Y,\kappa):=\projlim_{(\goth F,\succ)}(\Bbb R_+^{\Cal F},\nu^{\Cal F})$ is well defined as a standard probability space.
We now define a Borel map $\Phi:X^*\to Y$ by setting
$$
\Phi(x^*):=(\Phi(x^*)_{\Cal F})_{\Cal F\in\goth F},\quad\text{where } (\Phi(x^*)_{\Cal F}:=(N_K(x^*))_{K\in\Cal F}.
$$
Then $\Phi$ is one-to-one and onto.
Indeed, if two nonnegative Radon measures take the same values on every element of $\Cal K$ then these measures are equal.
This proves that $\Phi$ is one-to-one.
On the other hand, each element $y$ of $Y$ can be interpreted as a finitely additive non-negative measure on $\Cal K$, i.e. as a map $y:\Cal K\to\Bbb R_+$ such that $y(K_1\sqcup\dots\sqcup K_s)=y(K_1)+\cdots+y(K_s)$ for every collection  $\{K_1,\dots,K_s\}\in\goth F$.
Then $y$ extends uniquely to a $\sigma$-finite (non-negative) measure on $X$.
Of course, the extension is a Radon measure.
Hence $\Psi$ is onto.
Thus $\Psi$ is a Borel isomorphism of $X^*$ onto $Y$.
It remains to export $\kappa$ to $X^*$ via $\Phi^{-1}$ and
denote this measure by $\mu^*$.
It follows, in particular, that for each $K\in\Cal K$ and $j\in\Bbb Z_+$,
$$
\mu^*(\{x^*\in X^*\mid x^*(K)=j\})=\frac{\mu(K)^je^{-\mu(K)}}{j!}.\tag1-1
$$
If $X$ is partitioned into a union of mutually disjoint open non-compact subsets $X_1,\dots, X_m$ then
the mapping
$$
 x^*\mapsto(x^*\restriction X_1,\dots, x^*\restriction X_m)
 $$
is a measure preserving  Borel isomorphism of $(X^*,\mu^*)$ onto the direct product space $(X_1^*\times\cdots\times X_m^*,(\mu\restriction X_1)^*\times\cdots\times(\mu\restriction X_m)^*)$.

Given a Borel $\sigma$-algebra
$\goth F\subset\goth B$ generated by a topology which is weaker than the original topology on $X$, we denote by $\goth F^*$ the smallest sub-$\sigma$-algebra of $\goth B^*$
such that the mapping $N_K$ is measurable for each $K\in\goth F\cap\Cal K$.
For an  increasing sequence  of topologies $\tau_1\prec\tau_2\prec\cdots$ on $X$ which are weaker than the original topology on $X$, we denote by  $\goth F_1\subset\goth F_2\subset\cdots$ the corresponding increasing sequence of
Borel   sub-$\sigma$-algebras  generated by these topologies.
Since every compact open subset from $\bigvee_{n=1}^\infty\goth F_n$ is contained in $\goth F_m$ for some $m>0$, it folows that
 $$
 \left(\bigvee_{n=1}^\infty\goth F_n\right)^*=\bigvee_{n=1}^\infty\goth F^*_n.\tag1-2
 $$

Given an action $T=(T_g)_{g\in G}$ of $G$ on $X$ by $\mu$-preserving homeomorphisms $T_g$, we  associate a Borel action $T^*:=(T_g^*)_{g\in G}$ on $X^*$ by setting $T_g^*x^*:=x^*\circ T_g^{-1}$, $g\in G$.
It follows from \thetag{1-1} that $T^*$ preserves $\mu^*$.
The dynamical system $(X^*,\goth B^*,\mu^*, T^*)$ is called the {\it Poisson suspension} of $(X,\goth B,\mu,T)$.

\head 2. Proof of Theorem 0.1
\endhead

We will proceed in several  steps.

{\sl Step 1.}
First we construct a free strictly ergodic  infinite measure preserving $G$-action  on a locally compact non-compact Cantor space.
For that we will utilize the $(C,F)$-construction (see  \cite{Da1}, \cite{Da5}).
Let $(F_n,C_{n+1})_{n=0}^\infty$ be a sequence of finite subsets in $G$ such that
$(F_n)_{n=0}^\infty$ is a F{\o}lner sequence in $G$, $F_0=\{1\}$ and for each $n\ge 1$,
the following three basic conditions are satisfied:
\roster
\item"---"
$1\in F_n$ and $\# C_n>1$,
\item"---"
$F_n^{-1}F_nF_nC_{n+1}\subsetneq F_{n+1}$,
\item"---"
$F_nc\cap F_nc'=\emptyset$  for all $c\ne c'\in C_{n+1}$.
\endroster
We let $X_n:=F_n\times C_{n+1}\times C_{n+2}\times\cdots$.
Then $X_n$ endowed with the infinite product of the discrete topologies on $F_n$ and $C_j$, $j>n$, is a compact Cantor set.
Moreover, the mapping
$$
X_n\ni (f_n,c_{n+1}, c_{n+2},\dots)\mapsto(f_nc_{n+1}, c_{n+2},\dots)\in X_{n+1}
$$
is a topological embedding  of $X_n$ into $ X_{n+1}$.
We now consider the union $X=\bigcup_{n\ge 1}X_n$ and endow it with the topology of inductive limit.
Then $X$ is a locally compact non-compact Cantor set  and  $X_n$ is a compact open subset of $X$ for each $n$.
Given $g\in G$ and $n>0$, we  define a homeomorphism $T^{(n)}_g$ from a clopen subset
$(g^{-1}F_n\cap F_n)\times C_{n+1}\times C_{n+2}\times\cdots$ of $X_n$ to a clopen subset
$(F_n\cap gF_n)\times C_{n+1}\times C_{n+2}\times\cdots$ of $X_n$ by setting
$$
T^{(n)}_g(f_n,c_{n+1},\dots)\mapsto(gf_n, c_{n+1},\dots).
$$
It is easy to verify that the sequence $(T^{(n)}_g)_{n\ge 1}$ determines uniquely a homeomorphism $T_g$ of $X$ such that  $T_g\restriction X_n=T^{(n)}_g$ for each $n$.
It is straightforward to check that $T:=(T_g)_{g\in G}$ is a  free topological action of $G$ on $X$.
This action is minimal and uniquely ergodic, i.e. there is a unique $T$-invariant $\sigma$-finite Radon measure $\mu$ on $X$ such that $\mu(X_0)=1$.
To define $\mu$ explicitly we consider
for each $n\ge 0$  and $f\in F_n$,  a subset  $[f]_n:=\{(f,c_{n+1},c_{n+2},\dots)\in X_n\mid c_i\in C_i\text{ for all }i>n\}$.
Then $[f]_n$ is compact and open.
We call it a {\it cylinder}.
The family of all cylinders  is a base for the topology on $X$.
Every compact open subset of $X$ is a union of finitely many mutually disjoint cylinders.
 Hence every Radon measure on $X$ is determined uniquely by its values on the cylinders.
  It remains to note that  $\mu([f]_n)=\frac 1{\# C_1\cdots\# C_{n}}$
  for each $n\ge 0$ and $f\in F_n$.

In addition to the aforementioned three basic conditions on $(F_n,C_{n+1})_{n\ge 0}$,
we  will assume that 
$$
\lim_{n\to\infty}\frac{\# F_n}{\# C_1\cdots\# C_{n}}=\infty.
$$
This condition is equivalent to the fact that $\mu(X)=\infty$.

{\sl Step 2.} On this  step we obtain a free  probability preserving $G$-action.
For that we need two more conditions on  $(F_n, C_{n+1})_{n=0}^\infty$:
\roster
\item"$(\triangle)$" for each element $g$ of infinite order in $G$, there are infinitely many $n$ such that
$g^{l_n}F_nC_{n+1}\subset F_{n+1}\setminus(F_nC_{n+1})$ for some $l_n>0$.
\item"$(\square)$" for each element $g$ of finite order in $G$, there are infinitely many $n$ such that
$gF_n=F_n$.
\endroster
It is easy to see that if $(\triangle)$ is satisfied then 
$$ 
T_{g^{l_n}}X_n= T_{g^{l_n}}\bigsqcup_{f\in F_n}[f]_n=T_{g^{l_n}}\bigsqcup_{f\in F_n}\bigsqcup_{c\in C_{n+1}}[fc]_{n+1}=\bigsqcup_{f\in F_n}\bigsqcup_{c\in C_{n+1}}[g^{l_n}fc]_{n+1}.
$$
We used the fact that $X_n=\bigsqcup_{f\in F_n}[f]_n$ and $T_s[f]_n=[sf]_n$ whenever 
$f,sf\in F_n$ and $s\in G$.
Hence $T_{g^{l_n}}X_n\subset X_{n+1}$
 and $T_{g^{l_n}}X_n\cap X_n=\emptyset$.

Let $(X^*,\mu^*, T^*)$ be the Poisson suspension of $(X,\mu,T)$.
Then $\mu^*$ is $T^*$-invariant and $\mu^*(X^*)=1$.
We now verify that $T^*$ is free.
For that we will check that for each $g\in G\setminus\{1\}$, the subset of fixed points of the transformation $T^*_g$ is $\mu^*$-null.
Consider two cases.
If $g$ is of infinite order  and $A$ is a $T_g$-invariant subset of positive finite measure in $X$ then
in view of $(\triangle)$ there is $n>0$ such that $\mu(A\cap X_n)>0.9\mu(A)$ and $T_{g^{l_n}}X_n\cap X_n=\emptyset$.
Hence
$$
\mu(A\cap X_n)=\mu(T_{g^{l_n}}A\cap T_{g^{l_n}}X_n)=\mu(A\cap T_{g^{l_n}}X_n)\le\mu(A\setminus X_n)\le 0.1\mu(A),
$$
a contradiction.
Thus the transformation $T_g$ has no invariant subsets of finite positive measure.
Therefore the Poisson suspension $T_g^*$ of $T_g$ is weakly mixing \cite{Roy}.
This yields that $\mu^*(\{x^*\in X^*\mid T_g^*x^*=x^*\})=0$.
Consider now the case where $g$ is of finite order.
Let $H\subset G$ denote the cyclic subgroup generated by $g$.
It follows from  $(\square)$ that there is an open subset $Y\subset X$ of infinite measure such that
the sets $T_hY$, $h\in H$, form an open partition of $X$.
Indeed let $gF_{n_i}=F_{n_i}$ for an increasing sequence $n_1<n_2<\cdots$.
Choose a subset $S_1\subset F_{n_1}$ which meets each $H$-coset in $F_{n_1}$ exactly once. 
If $i>1$ choose a subset $S_i\subset F_{n_i}$ which meets each $H$-coset in $F_{n_i}\setminus(F_{n_{i-1}}C_{n_{i-1}+1}\cdots C_{n_i})$ exactly once.
It remains to set $Y:=\bigsqcup_{i=1}^\infty\bigsqcup_{s\in S_i}[s]_{n_i}$.
Of course, $Y$ is open and non-compact and $X=\bigsqcup_{h\in H}T_hY$.
Since $\infty=\mu(X)=\mu(\bigsqcup_{h\in H}T_hY)$ and $T$ preserves $\mu$, it follows that $\mu(Y)=\infty$.
Then the dynamical system $(X^*,\mu^*, (T_{h}^*)_{h\in H})$ is isomorphic to the finite direct product $(Y^*,(\mu\restriction Y)^*)^{H}$ endowed with the natural shiftwise action of $H$ (see \S 1).
Since the measure $(\mu\restriction Y)^*$ is nonatomic, this action is free.

{\sl Step 3.} 
We now verify that $h(T^*)=0$.
Let $\goth F_n$ denote the $\sigma$-algebra on $X$ generated by a single set $[1]_n$
 and let $\goth B_n$ denote the $\sigma$-algebra on $X$ generated by the compact open sets
 $[f]_n$, $f\in F_n$.
 Then $\goth B_1\subset\goth B_2\subset\cdots$ and $\bigvee_{n=1}^\infty\goth B_n$
is the entire Borel $\sigma$-algebra $\goth B$ on $X$.
Moreover,  $\goth B_n=\bigvee_{g\in F_n}T_g\goth F_n$ and
$$
\goth B_n^*=\bigvee_{g\in F_n}T_g^*\goth F_n^*\subset
 \bigvee_{g\in G}T_g^*\goth F_n^*=:(\goth F_n^*)^G
 $$
  for each $n$.
Therefore
$\bigvee_{n=1}^\infty(\goth F_n^*)^G\supset \bigvee_{n=1}^\infty\goth B_n^*=\goth B^*$.
The latter equality follows from \thetag{1-2}.
Moreover, it is easy to verify that  $(\goth F_1^*)^G\subset(\goth F_2^*)^G\subset\cdots$.
Hence
$$
h(T^*)=\lim_{n\to\infty} h\bigg(T^*\bigg|(\goth F_n^*)^G\bigg)\le\limsup_{n\to\infty}H(\goth F_n^*).
$$
We recall that $\Cal K$ denote the collection of all compact open subsets of $X$.
Since $\goth F_n\cap\Cal K=\{[1]_n\}$,
 the $\sigma$-algebra  $\goth F_n^*$ is generated by a countable partition of $X^*$ into the sets
 $N_{[1]_n}^{-1}(r)=\{x^*\in X^*\mid x^*([1]_n)=r\}$, $r=0,1,\dots$.
  This yields that
$$
\limsup_{n\to\infty}H(\goth F_n^*)=\lim_{n\to\infty} f(\mu([1]_n)),
$$
 where $f(t)$ is the entropy of the Poisson distribution  $(e^{-t}, e^{-t}t, e^{-t}t^2/2,\dots)$.
Since $\mu([1]_n)=1/(\#C_1\cdots\#C_{n})\to 0$,
it follows that $h(T^*)=0$.

{\sl Step 4.}
 On this step we show that some extra conditions on $(F_n,C_{n+1})_{n\ge 0}$ imply mixing for
the dynamical system $(X^*,\mu^*,T^*)$. 
Thus from now on we will assume that
 the following hold for each $n$ (in addition to  the conditions on $(F_n,C_{n+1})_{n\ge 0}$ listed above):
 \roster
\item"(i)"
$F_nF_n^{-1}F_nC_{n+1}\subset F_{n+1}$,
\item"(ii)"
the sets
$F_nc_1c_2^{-1}F_n^{-1}$, $c_1\ne c_2\in C_{n+1}$, and $F_nF_n^{-1}$ are all pairwise disjoint and
\item"(iii)"
$\# C_n\to\infty$ as $n\to\infty$.
\endroster
Denote by $U_T=(U_T(g))_{g\in G}$ the associated Koopman unitary representation of $G$ in $L^2(X,\mu)$, i.e. $U_T(g)f:=f\circ T_g^{-1}$ for each $f\in L^2(X,\mu)$.
As was shown in \cite{Da4, Theorem~5.1}\footnote{Though we considered in \cite{Da4} mainly the actions of the Heisenberg group $H_3(\Bbb Z)$, the proof  of Theorem~5.1 there does not use any specific property of $H_3(\Bbb Z)$. It holds for each amenable group.}, the conditions (i)--(iii) imply that $T$ is mixing as an infinite measure preserving action, i.e.
 $U_T(g)\to 0$ weakly as $g\to\infty$.
Let $U_{T^*}$ stand for the Koopman unitary representation of $G$ in $L^2(X^*,\mu^*)$ associated with the Poisson suspension $(X^*,\mu^*,T^*)$.
It is well known that $U_{T^*}$ is unitarily equivalent to the Fock unitary representation of $G$ generated by $U_T$
in the Fock space generated by $L^2(X,\mu)$ (see \cite{Ne}).
It follows that $U_{T^*}(g)$ converges weakly to the orthogonal projection to the constants in $L^2(X^*,\mu^*)$ as $g\to\infty$.
This is equivalent to the fact that $T^*$ is mixing.

Summarizing, we obtain that  $T^*$ is mixing free action of $G$ on $(X^*,\mu^*)$ 
and $h(T^*)=0$.
Thus Theorem~0.1 is proved.

\enddocument